\documentclass{llncs}
\usepackage{amssymb}
\usepackage{amsfonts}
\usepackage{latexsym}
\usepackage{graphics}
\def\Hide#1{\relax}
\newcounter{Seccounter}
\def\Section#1{{\bf \arabic{Seccounter}. #1.}%
	\addtocounter{Seccounter}{1}}
\newtheorem{Theorem}{Theorem}
\newtheorem{Lemma}{Lemma}
\newtheorem{Definition}{Definition}
\def\Proof{\par {\bf Proof: }}

\def\E#1{\mathbf{e}_{#1}}

\def\Fun#1#2{#1\rightarrow#2}

\def\ModII{\mbox{(\sf mod 2)}}
\def\Func#1{{\mathsf{#1}}}
\def\GA{\Func{GA}}

\def\Size#1{\# #1}

\def\N{\mathbb{N}}

\def\R{\mathbb{R}}
\def\Z{\mathbb{Z}}

\def\One{\mathbf{1}}

\def\x{\mathbf{x}}
\def\y{\mathbf{y}}
\def\z{\mathbf{z}}
\def\Fin{\mathbb{F}}
\def\symd{\ominus}
\def\diff{\mathop{\backslash}}
\def\|{\mathrel{|}}
\title{A Minimalist Construction of the Geometric Algebra}
\author{R.D. Arthan}
\institute{Lemma 1 Ltd. 2nd Floor, 31A Chain Street, Reading UK.  RG1 2HX\\
\email{rda@lemma-one.com}}
\begin{document}
\maketitle
\begin{abstract}
The geometric algebra is constructed from minimal raw materials.
\end{abstract}

\Section{Introduction}
Macdonald~\cite{Macdonald02} has noted that descriptions of the geometric algebra in the literature have tended either to use advanced concepts such as tensor algebra or to skip over or skimp on the proof of existence of the algebra.
He remedies this with a simple construction of the algebra based on a solution to the word problem for the defining equations of a Clifford algebra.
In this brief note, I give a construction along similar lines with a slightly different take on the solution to the word problem.
This simplifies Macdonald's approach in some respects, in particular, by using canonical representatives rather than equivalence classes.
As a technical convenience, I also adopt conventions that make the $\GA(p, q)$ for different $p$ and $q$ fit together so they can be constructed all at once.
As with Macdonald's construction, only an elementary knowledge of discrete mathematics and real vector spaces is required.

\Section{The Construction}
In the following, the variables $I$, $J$ and $K$ are always finite sets of integers,
$\Fin\Z$ denotes the set of all finite sets of integers, $\Size{A}$ denotes the number of elements of a finite set $A$, and $A \symd B$ denotes the symmetric difference, $(A \diff B) \cup (B \diff A)$, of any sets $A$ and $B$.

\begin{Definition}
Define functions $\alpha, \beta:\Fun{\Fin\Z \times \Fin\Z}{\N}$ and $\sigma:\Fun{\Fin\Z \times \Fin\Z}{\{+1, -1\}}$ as follows:
\begin{eqnarray*}
\alpha(I, J) &=& \Size{\{(i, j) \| i \in I \land j \in J \land j < i\}} \label{alpha}\\
\beta(I, J) &=& \Size{\{i \| i \in I \land i \in J \land i < 0\}} \label{beta}\\
\sigma(I, J) &=& (-1)^{\alpha(I, J) + \beta(I, J)} \label{sigma}
\end{eqnarray*}
\end{Definition}

\begin{Lemma} \label{lemmaI} For any $I, J, K \in \Fin\Z$, $\sigma(I, J)\sigma(I \symd J, K) = \sigma(I, J \symd K)\sigma(J, K)$
\end{Lemma}
\Proof
From the definition of $\sigma$, it suffices to prove the following congruences {\it modulo} 2:
\begin{eqnarray}
\alpha(I, J) + \alpha(I \symd J, K) &\equiv& \alpha(I, J \symd K) + \alpha(J, K) \ModII \label{congA}\\
\beta(I, J) + \beta(I \symd J, K) &\equiv& \beta(I, J \symd K) + \beta(J, K) \ModII
\label{congB}
\end{eqnarray}
From the definition of $\alpha$ and the inclusion/exclusion principle,
$\alpha(I \symd J, K) = \alpha(I, K) + \alpha(J, K) - 2\alpha(I \cap J, K)$.
So $\alpha(I \symd J, K)\equiv \alpha(I, K) + \alpha(J, K) \ModII$.
Similarly,
$\alpha(I, J \symd K) \equiv \alpha(I, J) + \alpha(I, K) \ModII$.
Adding these congruences and rearranging gives~(\ref{congA}).
Replacing $\alpha$ by $\beta$ in the same argument gives~(\ref{congB}).

\begin{Definition} Define $\GA(\infty, \infty)$ to be the real vector space with a basis of symbols $\E{I}$, one for each $I \in \Fin\Z$. Thus each element of $\GA(\infty, \infty)$ has a unique representation as a sum $\sum_{I \in \Fin\Z} \lambda_I\E{I}$, with real coefficients $\lambda_I$, where $\lambda_I = 0$ for all but finitely many $I$. Define a multiplication on $\GA(\infty, \infty)$ as follows:
\[
 (\sum_{I} \lambda_I\E{I}) (\sum_{J} \mu_J\E{J}) = \sum_{I, J} \lambda_I \mu_J \sigma(I, J) \E{I \symd J}
\]
Define $\One = \E{\{\}}$ and $\E{i} = \E{\{i\}}$ for each $i \in \Z$.
\end{Definition}

\begin{Theorem} \label{thmI}
$\GA(\infty, \infty)$ is an associative algebra over the field $\R$ of real numbers, having $\One$ as a two-sided unit. It is generated as an algebra by the elements $\E{i}$, which are linearly independent and satisfy:
\[
\begin{array}{lcrl}
\E{i}^2 &=& \One &\quad \mbox{if $i \ge 0$}\\
\E{i}^2 &=& -\One &\quad \mbox{if $i < 0$}\\
\E{i}\E{j} \hphantom{-} &=& -\E{j}\E{i} &\quad \mbox{if $i \not= j$}
\end{array}
\]
\end{Theorem}

\Proof
To say that $\GA(\infty, \infty)$ is an algebra over $\R$
is to say that it is a real vector space and that its multiplication satisfies the bilinearity conditions
$(\lambda\x)(\mu\y)=(\lambda\mu)(\x\y)$,
$(\x + \y)\z = \x\z + \y\z$ and
$\x(\y + \z) = \x\y + \x\z$
for any $\lambda, \mu \in \R$ and $\x, \y \in \GA(\infty, \infty)$.
It is a real vector space by construction and bilinearity is easily checked from the definitions.
The multiplication will be associative if $(\E{I}\E{J})\E{K} = \E{I}(\E{J}\E{K})$ for every $I, J, K \in \Fin\Z$,
i.e., if $\sigma(I, J)\sigma(I \symd J, K)\E{(I \symd J) \symd K} = \sigma(I, J \symd K)\sigma(J, K)\E{I \symd (J \symd K)}$,
but as $(I \symd J) \symd K = I \symd (J \symd K)$, this reduces to lemma~\ref{lemmaI}. It is easy to check that $\One\x=\x\One=\x$ for any $\x \in \GA(\infty, \infty)$, so $\One$ is indeed a two-sided unit for the algebra.
The set of elements of the form $\E{i}$ is a subset of the basis we have used to define $\GA(\infty, \infty)$ and so must be linearly independent.
By induction on $m$ and the definitions, one has that if
$I = \{i_1, i_2, \ldots, i_m\}$, with $i_1 < i_2 < \ldots < i_m$
then $\E{I} = \E{i_1}\E{i_2}\ldots\E{i_m}$. This means that the $\E{i}$ do indeed generate $\GA(\infty, \infty)$ as a real algebra.
Finally, the formulae for the various cases of $\E{i}^2$ and $\E{i}\E{j}$ are easy to verify from the definitions (for example, if $i < j$, then $\sigma(\{i\}, \{j\}) = 1$, and $\sigma(\{j\}, \{i\}) = -1$, so $\E{i}\E{j} = \E{\{i, j\}} = -\E{j}\E{i}$).

\begin{Definition} 
For $p, q \in \N$, define $\GA(p, q)$ to be the subalgebra of $\GA(\infty, \infty)$ generated by the $\E{i}$, $-q \le i < p$.
\end{Definition}

\begin{Theorem} \label{thmII}
$\GA(p, q)$ is an associative algebra over $\R$ with $\One$ as a two-sided unit. It has dimension $2^{p+q}$ as a real vector space and is generated as an algebra by elements $\E{i}$, $-q \le i < p$, which are linearly independent and satisfy:
\[
\begin{array}{lcrl}
\E{i}^2 &=& \One &\quad \mbox{if $i \ge 0$}\\
\E{i}^2 &=& -\One &\quad \mbox{if $i < 0$}\\
\E{i}\E{j} \hphantom{-} &=& -\E{j}\E{i} &\quad \mbox{if $i \not= j$}
\end{array}
\]
\end{Theorem}
\Proof
$\GA(p, q)$ has a basis comprising those $\E{I}$ where $I \subseteq \{-q, -q+1, \ldots, p-1\}$. As there are $2^{p+q}$ such $I$, $\GA(p, q)$ has dimension $2^{p+q}$. The rest is immediate from theorem~\ref{thmI}. 

\Section{Concluding Remarks}
No special properties of the field $\R$ have been used in this construction.
Any other field could have been used for the scalars.
Lundholm~\cite{Lundholm06} gives a generalised construction with some similarities to the one
presented here.

The equations in our two theorems define what is known as a Clifford algebra of type $(p, q)$.
Together with the condition on the dimension in theorem~\ref{thmII}, the equations characterise the $\GA(p, q)$ as what are called the universal Clifford algebras.
The real numbers themselves provide an example of a non-universal Clifford algebra of type $(1, 0)$: they satisfy the equations if one takes $\E{0} = \One = 1$, but in $\GA(1, 0)$, $\E{0} \not= \One$.
For more information see theorem 13.10 in \cite{Porteous69}.
The method of proof of that theorem shows that $\GA(\infty, \infty)$ is characterised by the conditions of theorem~\ref{thmI}.

The functions $\alpha$ and $\beta$ give the number of sign changes involved in putting a product of the elements $\E{i}$ in a normal form, see~\cite{Macdonald02}.
It is perhaps surprising that this intuition is scarcely relevant to the proofs given here.
The condition $i < j$ in the definition of the function $\alpha$ was only used above to derive the equations $\E{i}^2 = \pm1$ and $\E{i}\E{j} = -\E{j}\E{i}$.
Varying the condition gives similar algebras in which different pairs of basis elements can be made to commute or anticommute as desired.
However, none of these variants is a Clifford algebra.

It is traditional to give the $\E{i}$ that generate $\GA(p, q)$ positive indices, say with $\E{i}^2 = \One$ for $1 \le i \le p$ and  $\E{i}^2 = -\One$ for $p < i \le p + q$.
This has the disadvantage that $\GA(p, q) \not\subseteq \GA(p+r, q+s)$ unless $r=0$.
Our convention of using negative indices for the basis elements with negative squares means that the different $\GA(p, q)$ all fit together into the one large algebra $\GA(\infty, \infty)$.
$\GA(\infty, \infty)$ is perhaps a better candidate for the name ``{\em the} geometric algebra'' than any of the other $(\N \times \N)$-indexed family of contestants for the title. 

\bibliographystyle{plain}
\raggedright
\bibliography{bookspapers}

\end{document}